%% file: intr.tex
\documentclass[12pt]{article}
\usepackage{amstext,amsmath,amsfonts,graphics,color}
\newtheorem{thm}{Theorem} \newtheorem{lemma}[thm]{Lemma}
\newtheorem{prop}[thm]{Proposition} 
\newenvironment{pf} {\noindent{\sc Proof. }}{{\hfill
$\Box$}\par\vskip2\parsep} 
\newenvironment{pfof}[1]
{\par\vskip2\parsep\noindent{\sc Proof of\ #1. }}{{\hfill $\Box$}
  \par\vskip2\parsep}

\newcommand{\dof}{\bf\boldmath}

\newcommand{\lng}{\text{\rm long}}
\newcommand{\dist}{\text{\rm dist}}
\newcommand{\prob}{{\mbox{\bf P}\!}}
\newcommand{\zt}{{\mathbb Z}^2}
\newcounter{mycount}
\newenvironment{mylist}{\begin{list}{(\roman{mycount})}%
{\usecounter{mycount}\itemsep 0pt}}{\end{list}}

\title{Integrals, Partitions, and Cellular Automata}
\date{February 10, 2003 (revised May 4, 2003)}

\author{
Alexander E. Holroyd\thanks{Department of Mathematics,
University of British Columbia, Vancouver, BC, Canada V6T 1Z2, and 
Department of Mathematics, UC
Berkeley, CA 94720-3840, USA.  {\tt holroyd@math.berkeley.edu}. 
Research funded in part by NSF Grant DMS--0072398.}
\and Thomas M. Liggett\thanks{Department of Mathematics, UCLA, Los Angeles,
CA 90095-1555, USA. {\tt tml@math.ucla.edu}.  Research funded in part by NSF Grant 
DMS-00-70465.}
\and Dan Romik\thanks{Department of Mathematics,
                        Weizmann Institute of Science,
                        Rehovot 76100,
                        ISRAEL.
                        \tt{romik@wisdom.weizmann.ac.il}}}

\begin{document}

\maketitle
\renewcommand{\thefootnote}{}

\renewcommand{\thefootnote}{\arabic{footnote}} 

\begin{abstract}
We prove that 
$$\int_0^1\frac{-\log f(x)}xdx=\frac{\pi^2}{3ab}$$
where $f(x)$ is the decreasing function that satisfies $f^a-f^b=x^a-x^b$, for $0<a<b$.  When $a$ is an integer and $b=a+1$ we deduce several combinatorial results.  These include an asymptotic formula for the number of integer partitions not having $a$ consecutive parts, and a formula for the metastability thresholds of a class of threshold growth cellular automaton models related to bootstrap percolation.
\end{abstract}

\section{Introduction}
Let $0<a<b$ and define $f=f_{a,b}:[0,1]\rightarrow[0,1]$ to be
the decreasing function that satisfies
\begin{equation}
\label{deff}
[f(x)]^a-[f(x)]^b=x^a-x^b,\qquad 0\leq x\leq 1.
\end{equation}
Our central result is the following. 
\begin{thm}
\label{main}
 For every $0<a<b$,
$$\int_0^1\frac{-\log f(x)}xdx=\frac{\pi^2}{3ab}.$$
\end{thm}

\begin{figure}
\label{graph}
\caption{An illustration of the definition of the function $f$.}
\begin{center}
\input{graph.pstex_t}
\end{center}
\end{figure}
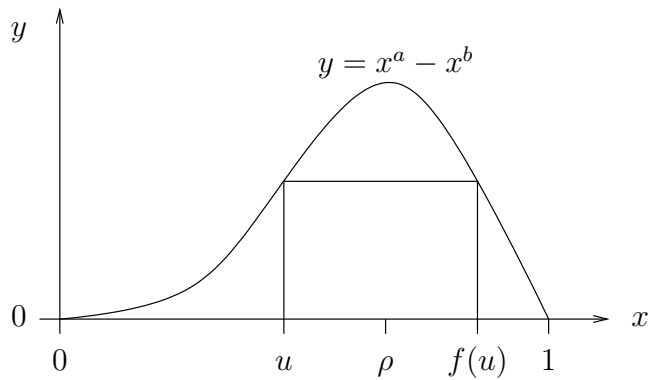

Note that $x^{a}-x^{b}$ is increasing on $[0,\rho]$
and decreasing on $[\rho,1]$ where $\rho=\rho_{a,b}\in(0,1)$ is defined by
\begin{equation}
\label{defrho}
-\log \rho =\frac{\log b -\log a}{b-a}
\end{equation}
It follows that $f$ is uniquely determined
by the above conditions, and satisfies
$$f(0)=1,\qquad f(1)=0,\qquad f(\rho)=\rho,$$
and 
$$f(f(x))=x,\qquad 0\leq x\leq 1.$$
See Figure \ref{graph}.

When $a$ is a positive integer and $b=a+1$, Theorem \ref{main} has the following consequences.

\paragraph{Probabilistic application.}
Let $0<s<1$, and let $C_1,C_2,\ldots$ be independent events with probabilities
$$\prob_s(C_n)=1-(1-s)^n$$
under a probability measure $\prob_s$.  (We can think of $C_n$ as the event that at least one occurs of a further set of $n$ independent events each of probability $s$).  Let $k$ be a positive integer, and let $A_k$ be the event
$$A_k=\bigcap_{i=1}^\infty (C_i \cup C_{i+1} \cup \cdots\cup C_{i+k-1})$$
that there is no sequence of $k$ consecutive $C_i$'s that do not occur.
\begin{thm}
\label{prob}
For every positive integer $k$,
$$-\log \prob_s(A_k) \sim \frac{\pi^2} {3k(k+1)}  \frac{1}{s} \qquad \text{as }s\to 0.$$
\end{thm}
The next two applications are consequences of Theorem \ref{prob}.

\paragraph{Number-theoretic application.}
A {\dof partition} of a positive integer $n$ is an unordered multiset of positive integers (called {\dof parts}) whose sum is $n$.  Let $p_k(n)$ be the number of partitions of $n$ that do not include any set of $k$ distinct consecutive parts.  (So for example $p_2(4)=4$, since the relevant partitions of $4$ are $(4),(3,1),(2,2),(1,1,1,1)$; the partition $(2,1,1)$ is not allowed because it has the $2$ consecutive parts $1,2$).
\begin{thm}
\label{part}
For every integer $k\geq 2$,
$$\log p_k(n) \sim \pi \sqrt{\frac{2}{3}\left(1-\frac{2}{k(k+1)}\right)n} \qquad \text{as }n\to \infty.$$
\end{thm}

\paragraph{Application to cellular automata.}
Threshold growth models are a class of simple cellular automaton models for nucleation and growth; see \cite{aizenman-leb},\cite{gravner-griffeath}, \cite{gravner-griffeath-2},\cite{h-boot} and the references therein.  
  Elements of the two-dimensional integer lattice $\zt$ are called {\dof sites}.  At each time step $t=0,1,2,\ldots $, a site is either {\dof active} or {\dof inactive}.  A site $z$ has a {\dof neighborhood} $N(z)\subseteq\zt$ defined by
$$N(z)=\{z+w:w\in N\},$$
where $N$($=N(0)$) is some fixed finite subset of $\zt$.  We also fix an integer $\theta$ called the {\dof threshold}.  The system evolves over time according to the following rules.
\begin{mylist}
\item A site that is active at time $t$ remains active at time $t+1$.
\item A site $z$ that is inactive at time $t$ becomes active at time $t+1$ if and only if its neighborhood $N(z)$ contains at least $\theta$ active sites at time $t$.
\end{mylist}

Consider the random initial state in which at time $0$, each site in the $L$ by $L$ square $\{1,\ldots ,L\}^2$ is active with probability $s$, independently for different sites, while all sites outside the square are inactive.  Let $I(L,s)$ be the probability that every site in the square eventually becomes active.  A central question is to determine for various models the behavior of the function $I(L,s)$ as $L\to\infty$ and $s\to 0$ simultaneously.
\begin{thm}
\label{boot}
Let $k\geq 2$ be an integer, and consider the threshold growth model with neighborhoods given by
$$N=N_k=\big\{(-v,0),(v,0),(0,-v),(0,v): v=1,2,\ldots,k-1\big\}$$
and threshold $\theta=k$.  For $L\to\infty$ and $s\to 0$ simultaneously we have
\begin{mylist}
\item
if $\liminf s\log L >\lambda$ then $I(L,s)\to 1$;
\item
if $\limsup s\log L <\lambda$ then $I(L,s)\to 0$,
\end{mylist}
where
$$\lambda=\frac{\pi^2}{3k(k+1)}.$$
\end{thm}

\paragraph{Further integrals.}
We can evaluate several other definite integrals using Theorem \ref{main}.
Recall the definition of $\rho$ in (\ref{defrho}).
\begin{thm}
\label{other1}
 For every $0<a<b$,
$$\int_0^\rho\frac{-\log f(x)}xdx=\frac{\pi^2}{6ab}-\frac{(\log \rho)^2}{2},$$
and
$$\int_\rho^1\frac{-\log f(x)}xdx=\frac{\pi^2}{6ab}+\frac{(\log \rho)^2}{2}.$$
\end{thm}
Define $\tilde{f}:[0,1]\rightarrow[0,1]$ to be
the decreasing function that satisfies
\begin{equation}
\label{defft}
-\tilde{f}(x)\log \tilde{f}(x)=-x\log x, \qquad 0\leq x\leq 1.
\end{equation}
\begin{thm}
\label{other2}
$$\int_0^1\frac{-\log \tilde{f}(x)}xdx=\frac{\pi^2}{3}.$$
\end{thm}
\begin{thm}
\label{other3}
$$\int_0^{e^{-1}}\frac{-\log \tilde{f}(x)}xdx=\frac{\pi^2}{6}-\frac{1}{2},$$
and
$$\int_{e^{-1}}^1\frac{-\log \tilde{f}(x)}xdx=\frac{\pi^2}{6}+\frac{1}{2}.$$
\end{thm}

\paragraph{Remarks.}
As we shall see in the proof of Theorem \ref{main}, the result for $f_{a,b}$ implies that for $f_{a\gamma,b\gamma}$ for any $\gamma>0$ via an easy argument.  The case $a=1,b=2$ is easy; in that case we have $f(x)=1-x$, and the integral in Theorem \ref{main} is standard (\cite{grad}, number 4.291.2).  The case $a=2,b=3$ also has an explicit formula for $f$, and this was used to prove Theorem \ref{main} in that case in \cite{h-boot}.  Our proof of the general case uses an entirely different approach.

The case $k=2$ of Theorem \ref{part} can also be deduced from a certain partition identity (see Section \ref{secpart}).  This raises the possibility of a family of partition identities corresponding to other values of $k$.  If such identities could be found, they might also lead to alternative (combinatorial) proofs of Theorems \ref{main},\ref{prob}.  We discuss these matters in more detail in Section \ref{secpart}.

The case $k=2$ of the threshold growth model in Theorem \ref{boot} is called {\em bootstrap percolation}.
Theorem \ref{boot} was proved for that case in \cite{h-boot}.  The general version is proved by a modification of the proof in \cite{h-boot}, making use of Theorem \ref{prob} above to obtain the numerical value of $\lambda$.  In Section \ref{secboot} we give an account of the proof, omitting some of the details.

Prior to the proof in \cite{h-boot}, even the existence of the sharp constant $\lambda$ in Theorem \ref{boot} was not known.  On the other hand, analogues of Theorem \ref{boot} with two {\em different} constants $\lambda_1,\lambda_2$ in (i),(ii) were known for a wide class of models.  In some cases, the ``scaling function'' $s\log L$ is replaced by a different function of $s,L$.  In particular, two-dimensional models are studied in detail in \cite{gravner-griffeath},\cite{gravner-griffeath-2}.  For neighborhoods as in Theorem \ref{boot} and threshold $k\leq \theta \leq 2k-2$ for example, the results in those articles imply that the appropriate scaling function is $s^{\theta-k+1} \log L$.  (The cases $\theta < k$ and $\theta > 2k-2$ turn out to be less interesting; in the former case, an active square of side $k$ will grow forever, while in the latter case an inactive square of side $k$ will remain inactive forever).  The reason for the particular choices of $N,\theta$ in Theorem \ref{boot} is that our methods (combined with those of \cite{h-boot}) yield the sharp constant $\lambda$ relatively easily in these cases.  The extension to other $N,\theta$ remains an open problem.

\section{Integrals}

In this section we prove Theorems \ref{main},\ref{other1},\ref{other2},\ref{other3}.

It suffices to prove the Theorem \ref{main} for the case $b=a+1$.
To check this, suppose that it holds for a
given choice of $a, b$.  For $\gamma>0$, let
$$g(x)=[f_{a,b}(x^{\gamma})]^{1/\gamma}.$$
Replacing $x$ with $x^{\gamma}$ in (\ref{deff}), we see that 
$$g^{a\gamma}(x)-g^{b\gamma}(x)=x^{a\gamma}-x^{b\gamma},$$
and $g$ is decreasing, so $g=f_{a\gamma,b\gamma}$.
Supposing that the theorem is true for $f_{a,b}$, we will 
check it for $g$:
\begin{eqnarray*}
\lefteqn{
\int_0^1\frac{-\log g(x)}{x}dx=   \frac 1{\gamma}\int_0^1
\frac{-\log f_{a,b}(x^{\gamma})}{x}dx}\\
&&=\frac
1{\gamma^2}\int_0^1\frac{-\log f_{a,b}(y)}{y}dy=\frac
1{\gamma^2}\frac{\pi^2}{3ab}=
\frac{\pi^2}{3(a\gamma )(b\gamma)}.
\end{eqnarray*}
In the second step above, we have made the 
change of variable $y=x^{\gamma}.$
So we may without loss of
generality take 
$$b=a+1.$$
Note that
in this case, $\rho=a/b$.
We will use
$\Gamma(\cdot)$ to denote the usual gamma function.

The proof of Theorem \ref{main} is based on properties of the
function
\begin{equation}
F(x)=\sum_{\ell=1}^{\infty}\frac{\Gamma(b\ell)}{\Gamma(a\ell)\ell!}
\frac{(x^a-x^b)^\ell}{a\ell}.\label{tml4}
\end{equation}
 By Stirling's
formula,
$$\frac{\Gamma(b\ell)}{\Gamma(a\ell)\ell!}\sim\bigg(\frac{b^{b}}{a^a}
\bigg)^\ell
\sqrt{\frac{a}{2\pi b \ell}},\qquad \ell\uparrow\infty.$$
Since the maximum value of $x^a-x^b=x^a(1-x)$ on $[0,1]$ is
$$\frac{a^a}{b^{b}},$$
the series in (\ref{tml4}) converges uniformly on $[0,1]$, and hence defines
a continuous function there. Note that the same cannot be
said for the series for $F'$. In fact, 
$F'$ is not continuous at $x=\rho$; one can
show using the proposition below that
$F'(\rho-)=1/\rho$ and $F'(\rho+)=-1/\rho$. This singularity
will play an important role in the analysis. The following
result contains the main properties of $F$ that will be
needed in the proof of Theorem \ref{main}.
\begin{prop} 
\label{ff}
Let $b=a+1$.  The function $F$ has
the following properties.
\begin{mylist}
\item $F(f(x))=F(x)$ on $[0,1]$.
\item $F(x)=-\log f(x)$ on $[0,\rho]$.
\item $F(x)=-\log x$ on $[\rho,1]$.
\item
$$\int_0^1\frac{F(x)}xdx=\frac{\pi^2}{6ab}.$$
\end{mylist}
\end{prop}

\begin{pf}
Part (i) is immediate from (\ref{deff}) and the fact that the
series in (\ref{tml4}) depends on $x$ through the expression
$x^a-x^{b}$. Part (ii) is a consequence
of (i) and (iii). To see this, take $x\in [0,\rho]$. Then
 $f(x)\in[\rho,1]$. By (iii), $F(f(x))=-\log f(x).$
Now use (i).  

Turning to the proof of (iii), define 
\begin{equation}
F(z)=\sum_{\ell=1}^{\infty}\frac{\Gamma(b\ell)}{\Gamma(a\ell)\ell!}
\frac{(z^a-z^b)^\ell}{a\ell}\label{tml5}
\end{equation}
for complex $z$ in the connected component $Q$ of  the
set
$$\bigg\{z\in {\mathbb C}: |z|^a|1-z|<\frac{a^a}{b^{b}}
\bigg\}$$
that contains the segment $(\rho,1]$. Note that
$Q$ is contained in the right half plane since for
$\mbox{Re}(z)=\rho$,
$|z|^a|1-z|\geq \rho^a(1-\rho)=a^a/b^b$. Therefore
$z^a$ can be defined unambiguously on $Q$ as an analytic
function that takes the value
1 at $z=1$, and $F$ is then analytic in $Q$. The
function $-\log z$ is also analytic in $Q$, and can
be chosen to take the value 0 at $z=1$. So, it suffices
to show that $F(z)=-\log z$ in a complex neighborhood of $z=1$.
Write $w=1-z$, and consider the neighborhood of 0
$$N=\bigg\{w\in {\mathbb C}:(1+|w|)^a|w|<\frac{a^a}{b^b}\bigg\}.$$
In $N$,
$$-\log(1-w)=\sum_{m=1}^{\infty}\frac{w^m}m.$$
Also, in $N$, the following rearrangement is justified
by absolute convergence of the series involved:
\begin{eqnarray*}
F(1-w)&=&\sum_{\ell=1}^{\infty}\frac{\Gamma(b\ell)}{\Gamma(a\ell)\ell!}
\frac{((1-w)^aw)^\ell}{a\ell}\\&=&\sum_{\ell=1}^{\infty}
\frac{\Gamma(b\ell)}{\Gamma(a\ell)\ell!}
\frac
1{a\ell}\sum_{k=0}^{\infty}\frac{\Gamma(a\ell+1)}{\Gamma(a\ell-k+1)k!}
(-1)^kw^{k+\ell}\\
&=&\sum_{m=1}^{\infty}b_mw^m,
\end{eqnarray*}
where 
$$b_m=\sum_{\ell=1}^m\frac{\Gamma(b\ell)(-1)^{m-\ell}}
{\ell!(m-\ell)!\Gamma(b\ell-m+1)}.$$
So, it suffices to prove that $b_m=1/m$ for $m\geq 1$.

To do so, use the property $\Gamma(\alpha+1)
=\alpha\Gamma(\alpha)$ to rewrite $b_m$ as
$$b_m=\sum_{\ell=1}^m\frac{(-1)^{m-\ell}}{\ell!(m-\ell)!}\prod_{i=1}^{m-1}
[b\ell-m+i].$$
The summand above that would correspond to $\ell=0$ is
$-1/m$. 
Therefore, $b_m=1/m$ is equivalent to
\begin{equation}
\sum_{\ell=0}^m\frac{(-1)^{m-\ell}}{\ell!(m-\ell)!}\prod_{i=1}^{m-1}
[b\ell-m+i]=0.\label{tml6}
\end{equation}
Now write
$$\frac{(1-x^b)^m}x=\sum_{\ell=0}^m\binom m\ell(-1)^\ell
x^{b\ell-1}.$$
To check (\ref{tml6}), it is then enough to show that
$$\frac{d^{m-1}}{dx^{m-1}}\frac{(1-x^b)^m}x\bigg|_{x=1}=0.$$
Let
$$h(x)=\frac 1x\bigg(\frac{1-x^b}{1-x}\bigg)^m,\qquad x\neq 1$$
and $h(1)=b^m$, so that
$$\frac{(1-x^b)^m}x=h(x)(1-x)^m.$$
Since $h$ is $C^{\infty}$ in a neighborhood of $x=1$,
$$\frac{d^{m-1}}{dx^{m-1}}h(x)(1-x)^m\bigg|_{x=1}=0$$
as required.

To prove part (iv) of the proposition, we use the standard beta
integral
$$\int_0^1u^{\alpha_1-1}(1-u)^{\alpha_2-1}du=\frac
{\Gamma(\alpha_1)\Gamma(\alpha_2)}{\Gamma(\alpha_1+\alpha_2)}$$
(see \cite{hps1} p148).
This gives
\begin{eqnarray*}
\int_0^1\frac{F(x)}xdx&=&\sum_{\ell=1}^{\infty}
\frac{\Gamma(b\ell)}{\Gamma(a\ell)\ell!}\int_0^1\frac{\big[x^a(1-x)\big]^\ell}
{a\ell}\frac
1xdx\\&=&\sum_{\ell=1}^{\infty}\frac{\Gamma(b\ell)\Gamma(a\ell)\Gamma(\ell+1)}
{\Gamma(a\ell)\ell!(a\ell)\Gamma(b\ell+1)}
\\&=&\frac1{ab}\sum_{\ell=1}^{\infty}\frac 1{\ell^2}=\frac{\pi^2}
{6ab}.
\end{eqnarray*}
\end{pf}

\begin{pfof}{Theorem \ref{main}}
As remarked at the beginning of the section, we may assume $b=a+1$.
 By Proposition \ref{ff} (iv), it
suffices to show that
\begin{equation}
\int_0^1\frac{-\log f(x)}xdx=2\int_0^1\frac{F(x)}xdx.\label{tml7}
\end{equation}
To do so, let
$$
I_1=\int_0^{\rho}\frac{-\log f(x)}xdx,\quad
I_2=\int_{\rho}^1\frac{-\log f(x)}xdx,$$
$$J_1=\int_0^{\rho}\frac{F(x)}xdx,\quad
J_2=\int_{\rho}^1 \frac{F(x)}xdx.$$
By Proposition \ref{ff} (ii), 
\begin{equation}
I_1=J_1.\label{tml8}
\end{equation}
Making the substitution $y=f(x)$ and then integrating by
parts gives
\begin{equation}
I_1=\int_{\rho}^1\log y\frac{f'(y)}{f(y)}dy
=-(\log\rho)^2+I_2,\label{tml9}
\end{equation}
since the boundary term at $y=1$ vanishes as a result of
$$f(y)\sim (1-y)^{1/a},\qquad y\uparrow 1.$$
By Proposition \ref{ff} (iii) and another integration
by parts,
\begin{equation}
J_2=\int_{\rho}^1\frac{-\log x}xdx=\frac {(\log\rho
)^2}{2}.\label{tml10}
\end{equation}
Combining (\ref{tml9}) and (\ref{tml10}) gives $I_2-I_1=2J_2$, and this,
together with (\ref{tml8}), gives $I_1+I_2=2J_1+2J_2$, which is (\ref{tml7}).
\end{pfof}

\begin{pfof}{Theorem \ref{other1}}
This follows immediately from (\ref{tml9}) above (which holds for all $0<a<b$) and Theorem \ref{main}.
\end{pfof}

\begin{pfof}{Theorems \ref{other2},\ref{other3}}
We will take the limit $a/b\to 1 $ in Theorems \ref{main},\ref{other1}.
Let $\epsilon\in(0,1)$, and write $f(x)=f_{1,1+\epsilon}(x)$ and 
$$\varphi(x)=\frac{x-x^{1+\epsilon}}{\epsilon},$$
so that by (\ref{deff}) we have 
\begin{equation}
\label{phif}
\varphi(f(x))=\varphi(x).
\end{equation}
Recall that $\varphi(x)$ is increasing on $[0,\rho]$ and decreasing on $[\rho,1]$.
Note that
$$\varphi(x)\to -x\log x \qquad\text{as }\epsilon\downarrow 0$$
uniformly in $x\in [0,1]$, and hence
$$f(x)\to \tilde{f}(x) \qquad\text{as }\epsilon\downarrow 0.$$
Note also (from (\ref{defrho})) that $\rho \downarrow e^{-1}$ as $\epsilon \downarrow 0$.

We will use dominated convergence.  First observe (by differentiating) that $\varphi(x)$ is decreasing in $\epsilon$ for each $x$.  Hence
$$x-x^2\leq \varphi(x) \leq -x\log x.$$
Therefore there exists a fixed constant $c$ satisfying
$$0<c<e^{-1}<\rho<1/2<1-c<1$$
such that for all $\epsilon \in(0,1)$ we have
$$u/2 \leq \varphi(u) \leq \surd u \qquad\text{for } u\in[0,c],$$
and
$$(1-u)/2 \leq \varphi(u) \leq 1- u \qquad\text{for } u\in[1-c,1].$$
It follows from (\ref{phif}) and the definition of $f$ that
there exist fixed positive constants $c',c''$ with $c'<1/4$ such that 
for all $\epsilon \in(0,1)$ we have
$$f(x)\geq\left\{
\begin{array}{ll}
1-2\surd x, \qquad & x<c' \\
c'', & c'\leq x \leq 1-c' \\
\left(\frac{1-x}{2}\right)^2, & x>1-c'.
\end{array}
\right. 
$$
Here, the bounds in the first and third cases are obtained by using the bounds for $\varphi$ above, and solving $(1-f(x))/2 \leq \surd x$ and $(1-x)/2 \leq \sqrt{f(x)}$ respectively.
Therefore for all $\epsilon \in(0,1)$ we have 
$$\frac{-\log f(x)}{x}\leq\left\{
\begin{array}{ll}
\frac{-\log (1-2\surd x)}{x}, \qquad & x<c' \\
c''', & c'\leq x \leq 1-c' \\
\frac{-2\log((1-x)/2)}{x} & x>1-c',
\end{array}
\right. 
$$
where $c'''>0$.
The function on the right is integrable on $[0,1]$ (see \cite{grad} for example).  Hence taking $\epsilon\downarrow 0$ and using the dominated convergence theorem,
Theorem \ref{other2} follows from Theorem \ref{main}, and 
Theorem \ref{other3} follows from Theorem \ref{other1}.
\end{pfof}

One may ask what Theorem \ref{main} yields in the limit when $a/b\to 0$ (or $\infty$).  In fact, it yields nothing new.  Taking $a\to 0$ with $ab=1$ say, an argument similar to the above shows that the limit of the integral is $\pi^2/3$, using only the (easy) case $a=1,b=2$ of Theorem \ref{main}.

\section{Probability}

In this section we prove Theorem \ref{prob}.

We say that a (finite or infinite) sequence of events has a {\dof $k$-gap} if there are $k$ consecutive events none of which occur.  Thus $A_k$ is the event that the sequence $C_1,C_2,\ldots $ has no $k$-gaps.

\begin{lemma}
\label{kgap}
Let $W_1,\ldots ,W_n$ be independent events each of probability $u\in(0,1)$.  Then the probability $g_n(u)$ that the sequence $W_1,\ldots ,W_n$ has no $k$-gaps satisfies
$$[f_{k,k+1}(1-u)]^n \leq g_n(u) \leq [f_{k,k+1}(1-u)]^{n-k+1}.$$
\end{lemma}

\begin{pf}
Writing $f=f_{k,k+1}(x)$, we have from (\ref{deff}) that
$$f^k - f^{k+1} = x^k - x^{k+1},$$
and rearranging gives
$$(f-x)f^k=(1-x)(f^k-x^k).$$
Provided $x\neq\rho$ we have $f\neq x$, so we may divide through by $f-x$ to obtain
\begin{equation}
\label{frr}
f^k=(1-x)(f^{k-1}+xf^{k-2}+x^2f^{k-3}+\cdots+x^{k-1}),
\end{equation}
and (\ref{frr}) holds when $x=\rho$ also by continuity (or (\ref{defrho})).

We now prove the statement of the lemma by induction on $n$.  For $n=0,\ldots,k-1$ we have $g_n(u)=1$, so the statement holds because $f_{k,k+1}(1-u)\in(0,1)$.  For $n\geq 0$, we may compute $g_{n+k}$ by conditioning on the first $W_i$ to occur:
\begin{eqnarray*}
g_{n+k}\!
&=&\!ug_{n+k-1}+(1-u)ug_{n+k-2}+(1-u)^2ug_{n+k-3}+\cdots+(1-u)^{k-1}ug_{n}\\
       &=&\!(1-x)(g_{n+k-1}+xg_{n+k-2}+x^2g_{n+k-3}+\cdots+x^{k-1}g_n),
\end{eqnarray*}
where we have written $x=1-u$.  Comparing this with (\ref{frr}) we deduce that if the lemma holds for $g_n,\ldots,g_{n+k-1}$ then it holds for $g_{n+k}$.
\end{pf}

\begin{pfof}{Theorem \ref{prob}}
The idea of the proof is that when $s$ is small, $\prob_s(C_n)$ varies only slowly with $n$, so we may use Lemma \ref{kgap} to deduce that $\prob_s(A_k)$ behaves approximately like $\prod_{n=1}^{\infty} f_{k,k+1}(1-\prob_s(C_n))$, and this in turn may be approximated using the integral in Theorem \ref{main}
(after a change of variable).

It is convenient to write
$$q=-\log(1-s)$$
so that $\prob_s(C_n)=1-e^{-nq}$ and $q\sim s$ as $s\to 0$.  Note that the indicator of $A_k$ is an increasing function of the indicators of $C_1,C_2,\ldots$, so if we increase (respectively decrease) the probabilities of the $C_i$ while retaining independence then we increase (respectively decrease) the probability of $A_k$.  We write 
$$r=\lfloor s^{-1/2} \rfloor$$
and let $C_n^+,C_n^-$ be independent events with probabilities 
\begin{eqnarray*}
\prob_s(C^+_n)&=&1-e^{-irq} \qquad (i-1)r < n\leq ir ;\\
\prob_s(C^-_n)&=&\left\{ 
 \begin{array}{ll}
s&\qquad 0<n\leq r, \\
1-e^{-irq}&\qquad ir<n\leq (i+1)r,
 \end{array}
\right.
\end{eqnarray*}
for $i=1,2,\ldots$.
Then  we have $\prob_s(C^-_n)\leq \prob_s(C_n)\leq \prob_s(C^+_n)$, and so $\prob_s(A^-_k)\leq \prob_s(A_k)\leq \prob_s(A^+_k)$, where $A^+_k$ (respectively $A^-_k$) is the event that the sequence $C^+_1,C^+_2,\ldots$ (respectively $C^-_1,C^-_2,\ldots$) has no $k$-gaps.

Now we may bound $\prob_s(A^+_k)$ above by the probability that
$$\text{for every $i\geq 1$, $C^+_{(i-1)r+1},\ldots,C^+_{ir}$ has no $k$-gaps.}$$
And we may bound $\prob_s(A^-_k)$ below by the probability that
\begin{eqnarray*}
&\text{$C^-_1,\ldots,C^-_r$ all occur, and}& \\
&\text{for every $i\geq 1$,
$C^-_{ir+1},\ldots,C^-_{ir+r-1}$ has no $k$-gaps, and $C^-_{ir+r}$ occurs}.&
\end{eqnarray*}
Hence, applying Lemma \ref{kgap} and writing $f=f_{k,k+1}$ we have
$$s^r\prod_{i=1}^\infty(1-e^{-irq})[f(e^{-irq})]^{r-1}\leq \prob_s(A_k)\leq \prod_{i=1}^\infty [f(e^{-irq})]^{r-k+1},$$
hence
\begin{eqnarray}
\lefteqn{(r-k+1)\sum_{i=1}^\infty -\log f(e^{-irq})\leq -\log\prob_s(A_k)} \nonumber \\
&&\leq   -r\log s +\sum_{i=1}^\infty -\log(1-e^{-irq}) +(r-1)\sum_{i=1}^\infty -\log f(e^{-irq}). \label{bounds1}
\end{eqnarray}

Applying the change of variable $x=e^{-z}$ to the integral in Theorem \ref{main} (with $a=k$, $b=k+1$) gives
$$\int_0^\infty -\log f(e^{-z}) dz=\frac{\pi^2}{3k(k+1)},$$
and in the special case $k=1$,
$$\int_0^\infty -\log (1-e^{-z}) dz=\frac{\pi^2}{6}.$$
Using the fact that $-\log f(e^{-z})$ and $-\log(1-e^{-z})$ are decreasing functions of $z$, (\ref{bounds1}) implies
\begin{eqnarray}
\lefteqn{\frac{r-k+1}{rq}\int_{rq}^\infty -\log f(e^{-z})dz\leq -\log\prob_s(A_k)} \label{bounds2} \\
&&\leq   -r\log s +\frac{1}{rq}\int_0^\infty -\log(1-e^{-z})dz +\frac{r-1}{rq}\int_0^\infty -\log f(e^{-z})dz. \nonumber
\end{eqnarray}
Now we let $s\to 0$.  Using the facts that $q\sim s$, $r\sim s^{-1/2}$, and both integrals are convergent, we obtain that the upper and lower bounds in (\ref{bounds2}) are both asymptotic to
$$\frac{1}{s}\int_0^\infty -\log f(e^{-z}) dz=\frac{\pi^2}{3k(k+1)}\frac{1}{s},$$
and hence the same holds for $-\log\prob_s(A_k)$.
\end{pfof}

\section{Partitions}
\label{secpart}

In this section we prove Theorem \ref{part}.

\begin{lemma}
\label{mono}
For any $k\ge 2$, $p_k(n)$ is a non-decreasing
function of $n$.
\end{lemma}

\begin{pf}
For $k\ge 3$, the following transformation defines
an injection of the set of partitions of $n$ not containing $k$
consecutive parts into the set of partitions of $n+1$ not containing
$k$ consecutive parts, thus establishing the claim. If the partition
does not contain all of the numbers $2,3,\ldots,k$ as parts, then we may
add another part equal to $1$, transforming the partition of $n$ into
a partition of $n+1$. If the partition \emph{does} contain $2,3,\ldots,k$
as parts, then we may transform it into a partition of $n+1$ by taking
one of the parts equal to $2$ and changing it into a $3$. It is easy
to verify that this is an injection.

It remains to prove the claim when $k=2$. For that case, we define the
following transformation taking partitions of $n$ without two
consecutive parts injectively into partitions of $n+1$ without two
consecutive parts.  If the partition does not contain any $2$'s, then
we may add a $1$. If the partition does contain $2$'s, we add $3$ to
the largest part in the partition and remove one $2$. (This fails for
the special partition $2=2$ of $n=2$, for that case verify the claim
directly). 
\end{pf}

\begin{pfof}{Theorem \ref{part}}
 Denote by 
$$G_k(x) = \sum_{n=0}^\infty p_k(n) x^n$$
the generating function of $p_k(n)$ ($k$ fixed). Let $p(n)$ be the
total number of (unrestricted) partitions of $n$, and denote its generating function
$$ G(x) = \sum_{n=0}^\infty p(n)x^n.$$
By \cite{newman}, p18, we have
\begin{equation}
\label{geng}
G(x) = \prod_{i=1}^\infty \frac{1}{1-x^i},\qquad 0<x<1.
\end{equation}
 We now observe
that $G_k(x)$ is closely related to the probability $\prob_s(A_k)$ in Theorem
\ref{prob}.  Let $s=1-x$.  We may
write the event $A_k$ as a disjoint union over 
the countable set ${\cal S}$ of all binary strings 
$a_1 a_2 a_3 a_4 \cdots \in \{0,1\}^\mathbb{N}$ that
contain only finitely many $0$'s, and in which
there are never $k$ consecutive $0$'s, of the event
$$\bigcap_{i\ :\ a_i = 1} C_i \cap \bigcap_{i\ :\ a_i = 0} C_i^{\text{c}}.$$
(By the Borel-Cantelli lemma, with probability one only finitely many
of the $C_i$'s will fail to occur). Therefore

\begin{eqnarray*}
\prob_s(A_k) 
&=& \sum_{a_1 a_2 a_3 \cdots \in {\cal S}}
 \prob_s\left(
\bigcap_{i\ :\ a_i = 1} C_i \cap \bigcap_{i\ :\ a_i = 0}
C_i^{\text{c}}\right) \\
& =& \sum_{a_1 a_2 a_3 \cdots \in {\cal S}}
\left[\prod_{i\ :\ a_i = 1} (1-x^i) \prod_{i\ :\ a_i = 0} x^i \right]
 \\
& = &
\left[\prod_{i=1}^\infty (1-x^i)\right]\sum_{a_1 a_2 a_3 \cdots \in {\cal S}}
\ \ \prod_{i\ :\ a_i=0}
 \frac{x^i}{1-x^i} \\
& = &
\frac{1}{G(x)} \sum_{a_1 a_2 a_3 \cdots \in{\cal S}}
\ \ \prod_{i\ :\ a_i=0}
 (x^i+x^{2i}+x^{3i}+\cdots ) \quad\text{(by (\ref{geng}))} \\
& = &\frac{G_k(x)}{G(x)},
\end{eqnarray*}
since on expanding out the sum of the products, the different ways to
 get $x^n$ correspond exactly to partitions of $n$ without $k$
 consecutive parts (choosing the power of $x^i$ corresponds to
 choosing the number of times the part $i$ appears in the partition).
Now using Theorem \ref{prob} and the
 standard fact (\cite{newman}, p19) 
$$ \log G(x) \sim \frac{\pi^2}{6(1-x)}, \qquad\text{as } x\uparrow 1,$$
we obtain $$ \log G_k(x) \sim
\frac{\pi^2}{6}\left(1-\frac{2}{k(k+1)}\right) \frac{1}{1-x}\qquad\text{as } x\uparrow 1.$$
We now use (a special case of) the Hardy-Ramanujan Tauberian Theorem \cite{hardy-ramanujan},
which says that if $H(x) = \sum_{n=0}^\infty b_n x^n$, where $b_n$ a
positive non-decreasing sequence, and $ \log H(x) \sim c/(1-x)$ as
$x\uparrow 1$, then $\log b_n \sim 2\sqrt{c n}$ as
$n\to\infty$. Theorem \ref{part} follows (using Lemma \ref{mono}).
\end{pfof}

\paragraph{The case $k=2$ and partition identities.} 
The special case $k=2$ of Theorem \ref{part} can be
deduced (and from it the corresponding cases of Theorems \ref{prob} and \ref{main})
using the following elementary partition identity due to P. A. MacMahon (\cite{andrews} p14,
examples 9, 10).

\begin{quote}
 The number of
partitions of $n$ not containing $1$'s and not containing two
consecutive parts is equal to the number of partitions of $n$ into
parts all of which are divisible by $2$ or $3$.
\end{quote}

Denote by $r(n)$ the number of such partitions of $n$. It is
straightforward to check that $r(n) \le r(n+2)$ for all $n$, since given a
no-ones, no-consecutive-parts partition of $n$ one may add 2 to its
largest part to turn it (injectively) into such a partition of $n+2$.  
Furthermore, we will
argue that the restriction of containing no 1's does not influence the
exponential rate of growth of the partition counting function, since we
have the inequalities 
\begin{equation} 
\label{same} 
\max\big\{r(n-1),r(n)\big\}\leq p_2(n) \leq \sum_{\ell=0}^n r(\ell). 
\end{equation} 
For the non-obvious part $r(n-1) \le p_2(n)$ of the lower bound, use the
(injective) transformation that adds 1 to the partition if there are no
2's, and otherwise takes a 2 and adds it, together with an additional 1,
to the largest part. For the upper bound, use the transformation that
takes a partition and deletes all the 1's.

Let
$R(x)=\sum_{n=0}^\infty r(n)x^n$ be the generating function of $r(n)$. By the
above partition identity we have for $0<x<1$
$$ R(x) = \prod_{i=0}^\infty \frac{1}{(1-x^{6i+2})(1-x^{6i+3})
(1-x^{6i+4})(1-x^{6i+6})}, $$
or
$$ \log R(x) = \sum_{j=2,3,4,6}\  \sum_{i=0}^\infty -\log(1-x^{6i+j}). $$
It can be shown in a manner analogous to the
asymptotic behavior of $G(x)$ cited above that for any $j\geq 1$,
$$ \sum_{i=0}^\infty -\log(1-x^{6i+j}) \sim \frac{1}{6}\cdot
\frac{\pi^2}{6(1-x)} \qquad\text{as }x\uparrow 1.$$
Therefore, summing over $j=2,3,4,6$ gives
$$ \log R(x) \sim \frac{4}{6}\cdot \frac{\pi^2}{6(1-x)}\qquad\text{as }x\uparrow 1.$$
By the Hardy-Ramanujan Tauberian Theorem applied to the function 
$R(x) + xR(x)$ (the generating function of the non-decreasing sequence 
$r(n-1)+r(n)$), we get
$$ \log \big(r(n-1)+r(n)\big) \sim \pi\frac{2}{3}\sqrt{n} 
\qquad\text{as }n\to\infty$$
which by (\ref{same}) gives Theorem \ref{part} for $k=2$.

 Now Theorem \ref{prob} may be deduced by following the
arguments of the proof of Theorem \ref{part} in the opposite direction (with
the slight adjustment of replacing $\prob_s(A_2)$ with
$\prob_s(C_1 \cap A_2)$ to account for the modified definition of the partitions),
and Theorem \ref{main} may be deduced by following the steps of the proof of
Theorem \ref{prob} (adapted to fit the modified statement) in the opposite
direction.  Note also that, as a consequence of MacMahon's identity, 
$\prob_s(C_1 \cap A_2)$ has the intriguing factorization
$$ \prob_s(C_1 \cap A_2) = \prod_{k=0}^\infty
(1-x^{6k+1})(1-x^{6k+5}) = (1-x)(1-x^5)(1-x^7)(1-x^{11}) \cdots $$ (where
again $x=1-s$). Can this fact be given a direct probabilistic proof?

We remark finally that, in light of the above argument and the neat
form of the exponential growth constant for $p_k(n)$ in Theorem
\ref{prob}, it is tempting to conjecture the existence of partition
identities for other integer values of $k$ that would give an
alternative proof of Theorem \ref{part} (and therefore of Theorems
\ref{prob} and \ref{main}) for integer $a=k$ and $b=a+1$. This would
imply, by analytic continuation, the general case of Theorem
\ref{main}, thus giving an independent proof of Theorem
\ref{main}. Presumably, such partition identities would equate the
number of partitions of $n$ not containing $k$ consecutive parts and
possibly satisfying some other ``mild'' conditions, with the number of
partitions of $n$ whose parts satisfy some congruence restrictions
modulo $k(k+1)$ (there should be two forbidden congruence classes) and
other mild restrictions. The discovery of such identities would be an
interesting positive use of partition asymptotics in the study of
partition identities. See \cite{brenner} for an example of a
\emph{negative} use of partition asymptotics, where they were used to
prove the non-existence of certain partition identities.
Also, see \cite{andrews2}, \cite{johansson}, \cite{johansson2} for
discussion of connections between partition theory and various models in
geometric probability, random matrix theory and statistical mechanics.

\section{Cellular Automata}
\label{secboot}

In this section we describe the proof of Theorem \ref{boot}.  The argument is a modification of that in \cite{h-boot}.  We therefore omit many of the details, concentrating instead on the differences compared with \cite{h-boot}.
The basic strategy is as follows.  Roughly, an $L$ by $L$ square will become fully active if and only if it contains at least one ``nucleation center''.  A nucleation center should be thought of as a local configuration of active sites that will grow to occupy the entire square.  One way for a nucleation center to occur involves the occurrence of (two independent copies of) the event $A_k$ in Theorem \ref{prob}.  Hence the estimate of $\prob_s(A_k)$ in Theorem \ref{prob} leads to the value of $\lambda$.  The above ideas lead to the bound Theorem \ref{boot} (i). 
Most of the work in \cite{h-boot} is in the proof of the bound (ii),
which involves ruling out the possibility of other types of nucleation center with substantially higher probabilities than the event described above.

The following set-up will be convenient.  Let $X$ be a random subset of $\zt$ in which each site is independently included with probability $s$.  More formally, denote by $\prob_s$ the product probability measure with parameter $s$ on the product $\sigma$-algebra of $\{0,1\}^{\zt}$, and define the random variable $X$ by $X(\omega)=\{x\in\zt: \omega(x)=1\}$ for $\omega\in\{0,1\}^{\zt}$.  A site $x\in\zt$ is said to be {\dof occupied} if $x\in X$.

Consider the threshold growth model as in Theorem \ref{boot}.
For a set of sites $K$, let $\langle K\rangle$ denote the set of eventually active sites if we start with $K$ as the set of initially active sites.  
We say that a set $K\subseteq\zt$ is {\dof internally spanned} if $\langle X\cap K\rangle=K$.
A {\dof rectangle} is a set of sites of the form
$R(a,b;c,d):=\{a,\ldots, c\}\times\{b,\ldots, d\}$; in particular
 we write $R(c,d)=R(1,1;c,d)$.  Thus 
$$I(L,s)=\prob_s\big(R(L,L) \text{ is internally spanned}\big).$$

\paragraph{Lower bound.}
Theorem \ref{boot} (i) is a consequence of the following lemma, which is the analogue of Theorem 2 (i) in \cite{h-boot}.
\begin{lemma}
\label{nuc}
$$\limsup_{s\to 0}\sup_{m\geq 1} -s\log I(m,s) \leq 2 \lambda.$$
\end{lemma}
(Here $\lambda$ is as in Theorem \ref{boot}).  Theorem \ref{boot} (i) is deduced from Lemma \ref{nuc} exactly as in the proof of Theorem 1 (i) in \cite{h-boot}, which in turn is based on an argument in \cite{aizenman-leb}.  The ideas are as follows; see \cite{h-boot} or \cite{aizenman-leb} for more details.  Let $S=S(L,s)$ be the event that $R(L,L)$ contains an internally spanned square of side $\lfloor s^{-3}\rfloor$.  The lemma implies easily that 
for $L,s$ as in Theorem \ref{boot} (i) we have $\prob_s(S)\to 1$.
On the other hand it may be shown that
$$\prob_s\big(R(L,L) \text{ is internally spanned} \mid S\big)\to 1,$$ 
proving the required result.
Next we turn to the proof of the lemma.
\begin{pfof}{Lemma \ref{nuc}}
Let $H=H(m)$ be the event that all the following occur:
$$\text{all sites in $R(k,k)$ are occupied},$$
$$\text{the sites $(m,1),(1,m)$ are occupied},$$
$$C_1,\ldots ,C_{m-k-1} \text{ has no $k$-gaps},$$
$$C'_1,\ldots ,C'_{m-k-1} \text{ has no $k$-gaps},$$
where $C_i$ is the event that $(k+i,j)$ is occupied for some $1\leq j\leq i$,
and $C'_i$ is the event that $(j,k+i)$ is occupied for some $1\leq j\leq i$.  See Figure \ref{nucleation}.
If $H$ occurs then $R(m,m)$ is internally spanned; to check this note that we may find an increasing sequence of internally spanned rectangles $R(i,j)$ with $i,j\in[k,m]$ and $|i-j|\leq k$, starting with $R(k,k)$ and ending with $R(m,m)$.  On the other hand, we have
$$\prob_s(H)\geq s^{k^2+2}\prob_s(A_k)^2$$
where $\prob_s(A_k)$ is as in Theorem \ref{prob}.  Therefore applying Theorem \ref{prob} yields the result.
\begin{figure}
\caption{An illustration of the event $H$ with $k=3$, $m=12$.}
\label{nucleation}
\begin{center}
\input{nuc.pstex_t}
\end{center}
\end{figure}
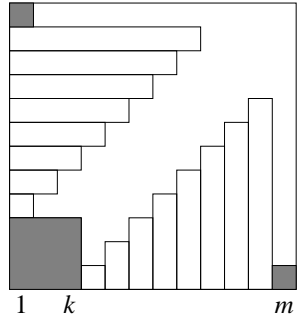
\end{pfof}

\paragraph{Upper bound.}
We now turn to the proof of Theorem \ref{boot} (ii).  
We start by defining a new cellular automaton model called the {\em enhanced model}.  We will explain how to prove the statement in Theorem \ref{boot} (ii) for the enhanced model, and this will imply the statement for the original model.

For a finite set of sites $F$, let ${\cal R}(F)$ be the smallest rectangle containing $F$.  Recall the definition of $N_k$ in Theorem \ref{boot}.  The {\dof enhanced model} is defined by the following rules.
\begin{mylist}
\item A site that is active at time $t$ remains active at time $t+1$.
\item For any site $z$, if $F$ is any set of $k$ sites in $N_k(z)$ that are active at time $t$, then all sites in ${\cal R}(F)$ are active at time $t+1$.
\item For any site $z'$, if $F'$ is any set of two sites in $N_2(z')$ that are active at time $t$, then $z'$ is active at time $t+1$.
\item Otherwise, inactive sites remain inactive.
\end{mylist}
We define $\langle \cdot\rangle$ and {\dof internally spanned} for the enhanced model in the same way as for the original model.
It is clear that any set that is internally spanned for the original model is internally spanned for the enhanced model.  Hence it is sufficient to prove that the statement in Theorem \ref{boot} (ii) holds for the enhanced model.  On the other hand, the enhanced model has the following useful property not shared by the original model:  if $K$ is a connected set of sites (in the sense of nearest neighbor connectivity on $\zt$) then $\langle K\rangle={\cal R}(K)$.  Furthermore, the rectangle ${\cal R}(F)$ in (ii) is of course connected.

We work with the enhanced model from now on, and our goal is to prove the statement in Theorem \ref{boot} (ii).
(As a consequence, both parts of Theorem \ref{boot} hold for both the enhanced model and the original model, with the same $\lambda$.  This is in contrast with the situation for critical probabilities of percolation models, where any ``essential enhancement'' of a model strictly lowers the critical probability; see \cite{aizenman-g}).
The proof follows the lines of \cite{h-boot}, and uses the following four lemmas, which we prove below.  The first three are deterministic properties of the enhanced model, and are analogues of Lemmas 9,10 and Proposition 30 in \cite{h-boot} respectively.

\begin{lemma}
\label{traversable}
For the enhanced model, if a rectangle is internally spanned then it has no $k$ consecutive unoccupied columns and no $k$ consecutive unoccupied rows.
\end{lemma}

The {\dof long side} of a rectangle $R=R(a,b;c,d)$ is defined as $\lng(R)=\max\{c-a+1,d-b+1\}$.
\begin{lemma}
\label{a-l}
Let $\ell$ be a positive integer and let $R$ be a rectangle satisfying $\lng(R)\geq \ell$. For the enhanced model, if $R$ is internally spanned then there exists an internally spanned rectangle $R'\subseteq R$ with $\lng(R')\in [\ell,2\ell+2k]$.
\end{lemma}

For rectangles $R,R'$ we write $\dist(R,R')=\min\{\|x-x'\|_\infty : x\in R,\; x'\in R' \}$.
\begin{lemma}
\label{disjoint}
Let $R$ be a rectangle with $|R|\geq 2$.  For the enhanced model, if $R$ is internally spanned then there exist distinct non-empty rectangles $R_1,\ldots ,R_r$, such that
\begin{mylist}
\item
$2\leq r\leq k$
\item
the {\em strict} inclusions $R_i\subset R$ hold for each $i$,
\item
$\langle R_1 \cup\cdots\cup R_r \rangle = {\cal R}(R_1 \cup\cdots\cup R_r) = R$,
\item
$\dist(R_i,R_j)\leq 2k$ for each $i,j$,
\item
$R_1,\ldots ,R_r$ are {\em disjointly} internally spanned; that is, there exist
disjoint sets of occupied sites $K_1,\ldots ,K_r$ with $\langle K_i\rangle=R_i$ for each $i$. 
\end{mylist}
\end{lemma}

\begin{lemma}
\label{concave}
For $k\geq 1$, the function $f_{k,k+1}$ is continuously differentiable and concave.
\end{lemma}

Once Lemmas \ref{traversable}--\ref{concave} are established, the proof of Theorem \ref{boot} (ii) for the enhanced model requires only minor modifications to the proof in \cite{h-boot}, so we omit the details.  The main modifications are as follows, referring to the terminology in that article.  The definition of a rectangle being ``horizontally (respectively vertically) traversable'' should be replaced with statement that no $k$ consecutive columns (respectively rows) are unoccupied.  The function ``$g(z)$'' in \cite{h-boot} should be replaced with $-\log f_{k,k+1}(e^{-z})$ in our notation (as in the proof of Theorem \ref{prob}).  The variational principles in Section 6 of \cite{h-boot} require that this function be convex; this follows from Lemma \ref{concave}.  Finally, in the definition of a ``hierarchy'', a vertex may have up to $k$ children (corresponding to having up to $k$ rectangles in Lemma \ref{disjoint} above), and a vertex is declared a ``splitter'' if it has two or more children.

Finally we give proofs of Lemmas \ref{traversable}--\ref{concave}. 
\begin{pfof}{Lemma \ref{traversable}}
Suppose that a rectangle $R$ has $k$ consecutive unoccupied columns, say.  Then it is easy to see that even if all other sites in $R$ is occupied, then no site in these columns lies in $\langle R\rangle$.
\end{pfof}

\begin{pfof}{Lemma \ref{a-l}}
This is a straightforward corollary of Lemma \ref{disjoint}.  Apply Lemma \ref{disjoint} to $R$, choose the resulting rectangle $R_i$ with the longest long side, then apply the lemma again to $R_i$, and so on, until the required $R'$ is obtained.
\end{pfof}

\begin{pfof}{Lemma \ref{disjoint}}
The proof is adapted from \cite{h-boot}.
Let $K=R\cap X$. By assumption we have $\langle K\rangle=R$.  We will 
construct $\langle K\rangle$ via the following algorithm.  We introduce a new time parameter $T$ not directly related to $t$ in the definition of the model.  For each time step $T=0,1,\ldots ,\tau$, we shall construct a collection of $m_T$ rectangles $R_1^T,\ldots , R_{m_T}^T$, and corresponding sets of sites $K_1^T,\ldots , K_{m_T}^T$, with the following properties:
\begin{mylist}
\item
$K_1^T,\ldots , K_{m_T}^T$ are pairwise disjoint;
\item
$K_i^T\subseteq K$;
\item
$R_i^T=\langle K_i^T\rangle$;
\item
if $i\neq j$ then $R_i^T \not\subseteq R_j^T$;
\item
$K\subseteq {\cal U}^T \subseteq \langle K\rangle$, where
$${\cal U}^T:=\bigcup_{i=1}^{m_T} R_i^T.$$
\end{mylist}

Initially, the rectangles and sets of sites are just the individual sites of $K$.  That is, let $K$ be enumerated as $K=\{x_1,\ldots, x_n\}$, and set 
$m_0=n$ \text{ and } $R_i^0=K_i^0=\{x_i\},$
so that in particular
$${\cal U}^0 =K.$$
The final set of rectangles will have the property that
\begin{equation}
\label{algfin}
{\cal U}^\tau =\langle K\rangle=R.
\end{equation}

We call a set of distinct rectangles $\{R_1,\ldots ,R_r\}$ a {\dof clique} if they satisfy conditions (i) and (iv) in Lemma \ref{disjoint}, and in addition $\langle R_1 \cup\cdots\cup R_r\rangle = {\cal R}(R_1 \cup\cdots\cup R_r)$.
The algorithm proceeds as follows.  Suppose $R_1^T,\ldots , R_{m_T}^T$ and $K_1^T,\ldots , K_{m_T}^T$ have already been constructed.
\begin{mylist}
\item[Step (I).]
If there does not exist a clique among $R_1^T,\ldots , R_{m_T}^T$, then {\bf stop}, and set $\tau=T$.

\item[Step (II).]  Suppose there does exist a clique; after reordering indices if necessary, denote it $R_1^T,\ldots R_r^T$.  Write $R'=\langle R_1^T \cup\cdots\cup R_r^T\rangle = {\cal R}(R_1 \cup\cdots\cup R_r)$, and $K'=K_1^T \cup\cdots\cup K_r^T$.
  
\item[Step (III).]  Construct the state $(R_1^{T+1},K_1^{T+1}),\ldots ,(R_{m_{T+1}}^{T+1},K_{m_{T+1}}^{T+1})$ at time $T+1$ as follows.  From the list $(R_1^{T},K_1^{T}),\ldots ,(R_{m_{T}}^{T},K_{m_{T}}^{T})$ at time $T$, delete every pair $(R_l^T,K_l^T)$ for which $R_l^T\subseteq R'$.  This includes $(R_i^T,K_i^T)$ for $i=1,\ldots,r$, and may include others.  Then add $(R',K')$ to the list.
\item[Step (IV).] Increase $T$ by 1 and return to Step (I).
\end{mylist}

It is straightforward to see that properties (i)--(v) are preserved by this procedure.  Also $m_T$ is strictly decreasing with $T$, so the algorithm must stop eventually.  We must check that (\ref{algfin}) is satisfied.  Suppose not; then there exists a site in $R\setminus{\cal U}^\tau$ which would become active in one step of the enhanced model if we start with ${\cal U}^\tau$ active.  Hence there exist $z,F$ or $z',F'$ as in rules (ii),(iii) of the enhanced model.  But the sites of $F$ or $F'$ must lie in at most $k$ different rectangles $R_1^\tau,\ldots ,R_r^\tau$ say; then it is easy to see that these rectangles form a clique, so the algorithm should not have stopped.
We claim that furthermore we must have $m_\tau=1$ and $R_1^\tau=R$.  If not, since ${\cal U}^\tau=R$, there must exist two distinct rectangles $R_i^\tau,R_j^\tau$ whose union is connected, but these two form a clique, so again the algorithm should not have stopped.

Finally, 
 considering the last time step of the algorithm (from time $\tau-1$ to time $\tau$) we obtain a set of rectangles with all the required properties.
\end{pfof}

\begin{pfof}{Lemma \ref{concave}}
It is sufficient to check the following
\begin{mylist}
\item
$f''(x)\leq 0 \text{ for }x\neq \rho,$
\item
$f'$ is continuous at $\rho$.
\end{mylist}
Let $\phi(x)=x^k(1-x)$, so that $\phi(x)=\phi[f(x)]$, and recall that 
$\phi'(x)=0$ only at $x=\rho=k/(k+1)$.
Differentiating $\phi(x)=\phi[f(x)]$ twice gives for $x\neq\rho$:
\begin{equation}
\label{deriv1}
f'(x)=\frac{\phi'(x)}{\phi'[f(x)]},
\end{equation}
\begin{eqnarray*}
 f''(x)&=&\frac{\phi''(x)-\phi''[f(x)][f'(x)]^2}{\phi'[f(x)]}
\\&
=&\frac{[\phi'(x)]^2}{\phi'[f(x)]}\big[\psi(x)-\psi[f(x)]
\big],
\end{eqnarray*}
where
$$\psi(x)=\frac{\phi''(x)}{[\phi'(x)]^2}=
\frac k{x^k}\frac{(k-1)-(k+1)x}{[k-(k+1)x]^2}.$$
We have $\phi'(f(x))<0$ for $x<\rho$ and $\phi'(f(x))>0$ for $x>\rho$,
therefore to prove (i) we need to show that
\begin{equation}
\label{tmlstar}
\psi(x)\geq\psi[f(x)]\quad\text{for}\quad 0\leq x<\rho
\end{equation}
and
$$\psi(x)\leq\psi[f(x)]\quad\text{for}\quad \rho< x\leq 1.$$
Since $f[f(x)]=x$, these two statements are equivalent, so we
will prove the first.

The first observation is that $\psi$ is decreasing on
$[0,\rho)$ and increasing on $(\rho,1]$. To see this,
compute
$$\psi'(x)=\frac{k}{x^{k+1}}\frac{k-1+(k+1)[(k+1)x-(k-1)]^2}
{[x(k+1)-k]^3},$$
and note that the numerator is nonnegative 
for all $x$ and the denominator
is negative for $x<\rho$ and positive for $x>\rho$. This immediately
gives (\ref{tmlstar}) for $0\leq x\leq (k-1)/(k+1)<\rho$:
$$\psi(x)\geq\psi\bigg(\frac{k-1}{k+1}\bigg)=0> -2k=\psi(1)
\geq \psi[f(x)].$$

To prove (\ref{tmlstar}) for $(k-1)/(k+1)\leq x<\rho$, we need the following
two facts:
\begin{equation}
\label{tml2star}
\phi(\rho-\sigma)\geq\phi(\rho+\sigma)\quad\text{and}
\quad\psi(\rho-\sigma)\geq\psi(\rho+\sigma)
\end{equation}
for $0<\sigma\leq 1/(k+1)$. Before checking this,
we will show that they imply 
(\ref{tmlstar}) for $(k-1)/(k+1)\leq x<\rho$.
Let $x=\rho-\sigma$, so that $0<\sigma\leq 1/(k+1)$. Then
$f(x)\leq\rho+\sigma=2\rho-x$ by the first statement in (\ref{tml2star}). Using
the second statement and the fact that $\psi$ is increasing on $(\rho,1]$ then gives
$$\psi(x)=\psi(\rho-\sigma)\geq\psi(\rho+\sigma)
\geq\psi[f(x)].$$

It remains to check (\ref{tml2star}). A little algebra shows that
$$\psi(\rho-\sigma)-\psi(\rho+\sigma)=
\frac{\rho}{\sigma^2}\frac{\phi(\rho-\sigma)-\phi(\rho+
\sigma)}{(\rho^2-\sigma^2)^k},$$
so the two statements in (\ref{tml2star}) are equivalent. To check the
first, compute
$$\frac d{d\sigma}\big[\phi(\rho-\sigma)-\phi(\rho+\sigma
)\big]=\sigma
(k+1)\big[(\rho+\sigma)^{k-1}-(\rho-\sigma)^{k-1}\big],$$
so that 
$\phi(\rho-\sigma)-\phi(\rho+\sigma)$ is increasing
in $\sigma$. Since this quantity is zero when $\sigma=0$,
it follows that it is nonnegative for $0<\sigma\leq 1/(k+1)$.

We have established (i).  This implies that the limits $f'(\rho+),f'(\rho-)$ exist (but are possibly infinite); to check (ii) it remains to show that they are equal.
Since $\phi''$ is continuous and non-zero at $\rho$, applying l'H\^{o}pital's rule to (\ref{deriv1}) gives
\begin{equation}
\label{square}
f'(\rho-)=\frac{\phi''(\rho-)}{\phi''(\rho+)f'(\rho-)}=\frac{1}{f'(\rho-)},
\end{equation}
And similarly for $f'(\rho+)$.  But $f$ is decreasing, so we must have $f'(\rho-)=f'(\rho+)=-1$.
\end{pfof}

\section*{Acknowledgements}

Alexander Holroyd thanks Laurent Bartholdi and Yuval Peres for valuable discussions.  We thank the referee for helpful comments.

\bibliography{phd}

\end{document}

%% file: graph.pstex_t
\begin{picture}(0,0)%
\includegraphics{graph.pstex}%
\end{picture}%
\setlength{\unitlength}{3947sp}%
\begingroup\makeatletter\ifx\SetFigFont\undefined%
\gdef\SetFigFont#1#2#3#4#5{%
  \reset@font\fontsize{#1}{#2pt}%
  \fontfamily{#3}\fontseries{#4}\fontshape{#5}%
  \selectfont}%
\fi\endgroup%
\begin{picture}(3904,2352)(879,-2930)
\put(1136,-2872){\makebox(0,0)[b]{\smash{\SetFigFont{12}{14.4}{\familydefault}{\mddefault}{\updefault}{\color[rgb]{0,0,0}$0$}%
}}}
\put(2543,-2872){\makebox(0,0)[b]{\smash{\SetFigFont{12}{14.4}{\familydefault}{\mddefault}{\updefault}{\color[rgb]{0,0,0}$u$}%
}}}
\put(3182,-2872){\makebox(0,0)[b]{\smash{\SetFigFont{12}{14.4}{\familydefault}{\mddefault}{\updefault}{\color[rgb]{0,0,0}$\rho$}%
}}}
\put(3759,-2872){\makebox(0,0)[b]{\smash{\SetFigFont{12}{14.4}{\familydefault}{\mddefault}{\updefault}{\color[rgb]{0,0,0}$f(u)$}%
}}}
\put(4206,-2872){\makebox(0,0)[b]{\smash{\SetFigFont{12}{14.4}{\familydefault}{\mddefault}{\updefault}{\color[rgb]{0,0,0}$1$}%
}}}
\put(4783,-2598){\makebox(0,0)[b]{\smash{\SetFigFont{12}{14.4}{\familydefault}{\mddefault}{\updefault}{\color[rgb]{0,0,0}$x$}%
}}}
\put(3246,-1000){\makebox(0,0)[b]{\smash{\SetFigFont{12}{14.4}{\familydefault}{\mddefault}{\updefault}{\color[rgb]{0,0,0}$y=x^a-x^b$}%
}}}
\put(879,-772){\makebox(0,0)[b]{\smash{\SetFigFont{12}{14.4}{\familydefault}{\mddefault}{\updefault}{\color[rgb]{0,0,0}$y$}%
}}}
\put(879,-2598){\makebox(0,0)[b]{\smash{\SetFigFont{12}{14.4}{\familydefault}{\mddefault}{\updefault}{\color[rgb]{0,0,0}$0$}%
}}}
\end{picture}

%% file: nuc.pstex_t
\begin{picture}(0,0)%
\includegraphics{nuc.pstex}%
\end{picture}%
\setlength{\unitlength}{1973sp}%
\begingroup\makeatletter\ifx\SetFigFont\undefined%
\gdef\SetFigFont#1#2#3#4#5{%
  \reset@font\fontsize{#1}{#2pt}%
  \fontfamily{#3}\fontseries{#4}\fontshape{#5}%
  \selectfont}%
\fi\endgroup%
\begin{picture}(3624,3912)(1189,-4261)
\end{picture}